\documentclass[11pt]{article}
\usepackage[utf8]{inputenc}
\usepackage{algpseudocode}
\usepackage{xcolor}
\usepackage{amsmath}
\usepackage{amssymb}
\usepackage{siunitx}
\usepackage{footnote}
\usepackage{float}
\usepackage{mathtools}
\usepackage{color,soul}
\usepackage[numbers,sort&compress]{natbib}
\usepackage[edges]{forest}
\usepackage{setspace}
\usepackage{algorithm2e}
\SetKwComment{Comment}{/* }{ */}
\RestyleAlgo{ruled}
\DontPrintSemicolon

\usepackage{doi}

\usepackage{authblk}
\usepackage[symbol]{footmisc}

\usepackage{orcidlink}

\usepackage[english]{babel}
\usepackage[font=small,labelfont=bf]{caption}
\usepackage[skip=8mm]{caption} 
\usepackage{subcaption}
\usepackage{graphicx}

\usepackage[short,nocomma,c1]{optidef}

\usepackage{duckuments}

\usepackage{tikz}
\usetikzlibrary{arrows,calc,fit,patterns}
\usepackage{eurosym}

\usepackage{pdfcomment}

\usepackage{enumerate}
\usepackage[shortlabels]{enumitem}

\usepackage{hyperref}
\usepackage{amsthm}
\usepackage[capitalize,nameinlink]{cleveref}
\usepackage{soul}
\usepackage{tabularx}

\usepackage[left=2cm,top=1.4cm,right=2cm,bottom=2.5cm,nohead]{geometry}
\usepackage{amssymb}
\usepackage{amsmath}
\usepackage{bbold}

\makeatletter
\renewcommand{\section}{\@startsection{section}{1}{0ex}%
                                   {-3.5ex \@plus -1ex \@minus -.2ex}%
                                   {0.1ex \@plus.2ex}%
                                   {\normalfont\Large\bfseries\sffamily}}
\renewcommand{\subsection}{\@startsection{subsection}{2}{0ex}%
                                     {-3.25ex\@plus -1ex \@minus -.2ex}%
                                     {1.5ex \@plus .2ex}%
                                     {\normalfont\large\bfseries\sffamily}}
\renewcommand{\subsubsection}{\@startsection{subsubsection}{3}{0ex}%
                                     {-3.25ex\@plus -1ex \@minus -.2ex}%
                                     {1.5ex \@plus .2ex}%
                                     {\normalfont\normalsize\bfseries\sffamily}}
\renewcommand{\paragraph}{\@startsection{paragraph}{4}{\z@}%
                                    {3.25ex \@plus1ex \@minus.2ex}%
                                    {-1em}%
                                    {\normalfont\normalsize\bfseries\sffamily}}
\renewcommand{\subparagraph}{\@startsection{subparagraph}{5}{\parindent}%
                                       {3.25ex \@plus1ex \@minus .2ex}%
                                       {-1em}%
                                      {\normalfont\normalsize\bfseries\sffamily}}
\renewcommand{\@maketitle}{%
  \newpage
  \null
  \begin{center}%
  \let \footnote \thanks
    {\Large \textsf{\textbf{\@title}} \par}%
    \vskip 0.5em%
    {\large
      \lineskip .5em%
      \begin{tabular}[t]{c}%
        \textsl{\@author}
      \end{tabular}\par}
  \end{center}%
  \par
  \vskip 1em}
\makeatother

\usepackage{tikz}

\title{When does learning pay off? A study on DRL-based dynamic algorithm configuration for carbon-aware scheduling}

\author[,a,b]{Andrea Mencaroni\,\orcidlink{0000-0002-0110-3218}\thanks{Corresponding author. \\ \textit{E-mail addresses}: 
\href{mailto:andrea.mencaroni@ugent.be}{andrea.mencaroni@ugent.be} (A. Mencaroni),
\href{mailto:r.v.j.reijnen@tue.nl}{r.v.j.reijnen@tue.nl} (R. Reijnen),
\href{mailto:yqzhang@tue.nl}{yqzhang@tue.nl} (Y. Zhang),
\href{mailto:dieter.claeys@ugent.be}{dieter.claeys@ugent.be} (D. Claeys).
}}
\author[c,d]{Robbert Reijnen\,\orcidlink{0000-0002-1629-6040}}
\author[c,d]{Yingqian Zhang\,\orcidlink{0000-0002-5073-0787}}
\author[a,b]{Dieter Claeys\,\orcidlink{0000-0002-7666-2479}}
\affil[a]{\footnotesize Department of Industrial Systems Engineering and Product Design, Ghent University, Ghent, Belgium} 
\affil[b]{\footnotesize Industrial Systems Engineering (ISyE), Flanders Make, Kortrijk, Belgium}
\affil[c]{\footnotesize Eindhoven Artificial Intelligence Systems Institute, Eindhoven University of Technology, The Netherlands}
\affil[d]{\footnotesize Department of Industrial Engineering \& Innovation Sciences, Eindhoven University of Technology, The Netherlands}
\date{} 

\providecommand{\keywords}[1]{\textbf{\textit{Keywords --}} #1}

\usepackage{xpatch}
\makeatletter
\AtBeginDocument{\xpatchcmd{\@thm}{\thm@headpunct{.}}{\thm@headpunct{}}{}{}}
\makeatother

\begin{document}

\maketitle

\section*{Abstract}
Deep reinforcement learning (DRL) has recently emerged as a promising tool for Dynamic Algorithm Configuration (DAC), enabling evolutionary algorithms to adapt their parameters online rather than relying on static tuned configurations.
While DRL can learn effective control policies, training is computationally expensive.
This cost may be justified if learned policies generalize, allowing a single training effort to transfer across instance types and problem scales.
Yet, for real-world combinatorial optimization problems, it remains unclear whether this promise holds in practice and under which conditions the investment in learning pays off.

In this work, we investigate this question in the context of the carbon-aware permutation flow-shop scheduling problem.
We develop a DRL-based DAC framework and train it exclusively on small, simple instances.
We then deploy the learned policy on both similar and substantially more complex unseen instances,
and compare its performance against a static tuned baseline, tuned under the same budget and on the same instances, which provides a fair point of comparison.

Our findings show that the proposed method provides a strong dynamic algorithm control policy that can be effectively transferred to different unseen problem instances. Notably, on simple and cheap to compute instances, similar to those observed during training and tuning, DRL performs comparably with the statically tuned baseline. However, as instance characteristics diverge and computational complexities increase, the DRL-learned policy continuously outperforms static tuning. These results confirm that DRL can acquire robust and generalizable control policies which are effective beyond the training instance distributions. This ability to generalize across instance types makes the initial computational investment worthwhile, particularly in settings where static tuning struggles to adapt to changing problem scenarios.

\keywords{Evolutionary algorithms, Deep reinforcement learning, Dynamic algorithm configuration, Carbon-aware scheduling}

\section{Introduction}
\label{section:introduction}

Evolutionary algorithms (EAs) have proven to be flexible and powerful approaches for solving complex combinatorial optimization problems \cite{Marti2025}.
Their performance, however, critically depends on the algorithm configuration, that is, the selection of variation operators and associated control parameters that govern the search dynamics \cite{Iomazzo2023}.
These elements jointly determine the balance between exploration and exploitation during the search process \cite{Eiben2011}, which is known to be crucial for EA performance \cite{Crepinsek2013}.
While greater exploration reduces the risk of premature convergence, stronger exploitation accelerates convergence around promising regions but increases the risk of missing better solutions elsewhere.

Traditionally, algorithm configuration is addressed through an offline tuning process, in which the algorithm is repeatedly executed on a set of training instances, and its performance is evaluated under different configurations.
The configuration yielding the best average performance within the available computational budget is then adopted for deployment on new instances of a similar type \cite{Eiben2012}.
Since such tuning processes are highly time-consuming \cite{Huang2020}, automated algorithm configuration methods have been developed to systematically explore the parameter space and identify well-performing configurations \cite{LopezIbanez2016, Akiba2019}.
Yet, even with automation, tuning remains computationally expensive and fundamentally instance-specific as configurations optimized on one instance distribution typically do not transfer reliably to others.
When the problem size or structure changes, tuning must be typically repeated from scratch to avoid significant performance degradation.

A second limitation of this paradigm is that it produces static configurations, which are inherently suboptimal because the optimal balance between exploration and exploitation evolves during the search \cite{Cervantes2009,Goldman2011}.
Early stages, in fact, typically benefit from broader exploration, whereas later stages require stronger exploitation to refine best-found solutions \cite{Li2010}.
This observation has motivated a shift towards dynamic or online parameter control, where parameter values and operators are adjusted during the execution of the algorithm \cite{Karafotias2015}.
While such methods could potentially be superior by letting the exploration-exploitation balance evolve dynamically, they are usually based on handcrafted rules or simple reactive mechanisms, which also remain problem- and instance-specific \cite{Karafotias2015}.
Consequently, as for static tuning, these methods are not designed or evaluated from a cross-instance generalization perspective, and typically require costly retuning or redesign whenever instance size or structure changes.

In recent years, learning-based approaches to online algorithm control have gained growing attention \cite{Adriaensen2022}.
These methods cast algorithm configuration as a sequential decision-making problem, where parameter values correspond to actions, and rewards incentivize changes that led to improved solutions \cite{Powell2022}.
Within this line of work, Dynamic Algorithm Configuration (DAC) has been introduced as a formal framework that models algorithm configuration as a Markov Decision Process (MDP), enabling the use of reinforcement learning to learn control policies \cite{Biedenkapp2020}.
In this way, an agent learns to adapt the algorithm's behavior online so as to optimize performance over the course of the search.

Beyond their potential to outperform static or rule-based strategies on individual instance types \cite{Chen2020, Quevedo2021, Cui2024}, the DAC framework offers the possibility of learning transferable control policies that generalize across problem instances and scales \cite{Biedenkapp2020}.
This stands in contrast to both static tuning and classical dynamic control strategies, which must be retuned or redesigned for each new instance type, and for which no form of cross-instance generalization has been demonstrated.
If such generalization can indeed be achieved, control policies could be trained on small instances, where training is computationally cheaper, and subsequently deployed on larger ones.
Under such a setting, training could be considered as a one-time investment justified by the potential for cross-instance generalization and the ability to transfer to more complex problem variants \cite{Sharma2019, Reijnen2024}.

At the same time, learning such policies requires substantial computational effort, often comparable to that of classical tuning approaches \cite{Reijnen2022}.
Moreover, as problem size increases (e.g., more machines in scheduling problems or more customers in routing problems), each algorithm run becomes longer \cite{He2001}, further increasing the computational burden of training.
Consequently, the practical usefulness of DAC critically depends on whether the learned policies exhibit sufficient generalization, so that training does not need to be repeated for every new instance type or scale.

Despite this promise, empirical evidence supporting cross-instance generalization remains limited.
While related lines of research, such as transfer learning and meta-learning in automated algorithm design \cite{Lu2015, Ma2026}, aim to facilitate the transfer of knowledge across instances, they primarily focus on developing general learning frameworks or on improving performance across tasks, rather than providing a systematic empirical assessment of the generalization of learned dynamic control policies.
Most existing studies evaluate their methods on relatively simple continuous optimization problems or standard benchmark functions \cite{Yang2025}, making it difficult to assess the benefits of these methods for real-world combinatorial optimization problems \cite{Reijnen2025}.
Furthermore, the majority of works focus on outperforming previous results on specific problem instances, i.e., playing the \textit{up-the-wall} game \cite{Burke2009}, without providing broader insight into when such methods are effective and when they are not.
To the best of our knowledge, only a limited number of studies have explicitly examined the generalization capabilities of learning-based DAC \cite{Reijnen2025}.
Yet these works do not investigate the conditions under which learning-based DAC is actually beneficial compared to classical parameter-tuning approaches.
Consequently, the extent to which learning-based DAC can generalize beyond its training distribution, and the conditions under which it is preferable to classic static configurations, remain unclear.

To investigate this potential, we develop a deep reinforcement learning (DRL)-based DAC framework that autonomously adapts parameters during execution.
We apply this framework to the carbon-aware scheduling (CAS) problem \cite{Mencaroni2025}, a challenging real-world optimization problem in which scalability is critical due to its time-dependent objective function.

Our main contributions are as follows:
\begin{itemize}
    \item We propose MA-DRL-CAS-PFSP, a framework that uses DRL to dynamically configure the parameters of a memetic algorithm (MA) for the permutation flow-shop problem (PFSP) under a carbon-aware scheduling objective.
    \item We design the training process and experiments to explicitly test whether and to what extent learned policies generalize and transfer to larger, real-world instance types.
    \item We compare the proposed dynamic approach against its static counterpart, thereby identifying the conditions under which using learning-based DAC is beneficial.
\end{itemize}

\section{Problem setting and related work}
\label{section:problem_setting_and_related_work}

Before reviewing the related work on DAC in evolutionary algorithms, we now summarize the key aspects of the CAS problem.

\subsection{Carbon-aware scheduling}
\label{subsection:CAS}
The CAS problem extends classical production scheduling problems by explicitly accounting for the carbon emissions associated with electricity consumption in manufacturing.
Unlike traditional scheduling approaches that typically focus on makespan minimization \cite{Ostermeier2024}, CAS introduces additional time-dependent factors such as individual job power requirements, grid carbon intensity, and on-site renewable electricity generation to minimize scope 2 greenhouse gas (GHG) emissions.
This time-dependency increases the computational complexity relative to classical scheduling problems \cite{Gawiejnowicz2025}, making CAS particularly sensitive to scalability challenges, especially when tuning is required for large-scale instances.
Consequently, CAS provides an ideal context for evaluating the scalability and effectiveness of adaptive parameter control in EAs.

While we refer to \citet{Mencaroni2025} for the full model details, we summarize the main characteristics here for ease of reference.
Operations require not only processing time but also electrical power to be executed on machines, and each operation is therefore associated with time-specific power requirements.
Machines are primarily powered by on-site renewable energy sources, which are preferred due to their zero scope 2 emissions, although their availability is limited and varies over time.
When renewable generation is insufficient, electricity is drawn from the public grid, which can be considered unlimited in availability but has a non-zero and time-varying carbon intensity depending on the grid's energy mix.
The objective is to minimize the scope 2 greenhouse gas (GHG) emissions while ensuring that all the jobs are completed within the given planning horizon.

\begin{figure}[h]
	\begin{subfigure}{.49\columnwidth}
		\centering
		\includegraphics[clip, trim={0.4cm 0cm 0.4cm 0.5cm}, width=\textwidth]{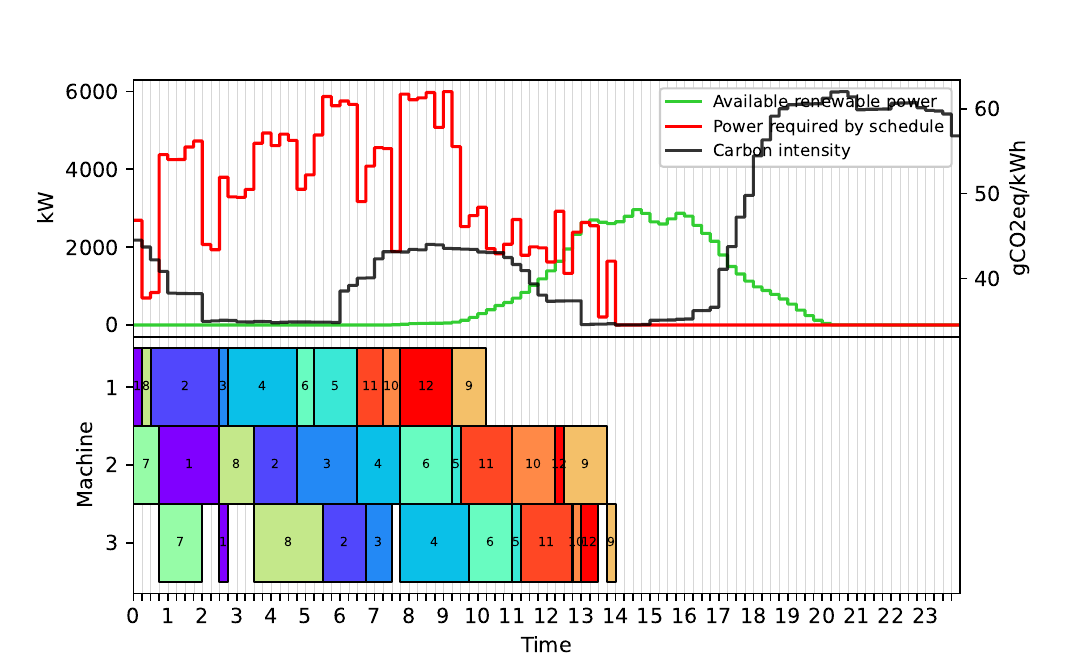}
		\caption{Minimization of makespan: $7.00$ tCO\textsubscript{2}}
		\label{fig:gantt_makespan}
	\end{subfigure}
	\begin{subfigure}{.49\columnwidth}
		\centering
		\includegraphics[clip, trim={0.4cm 0cm 0.4cm 0.5cm}, width=\textwidth]{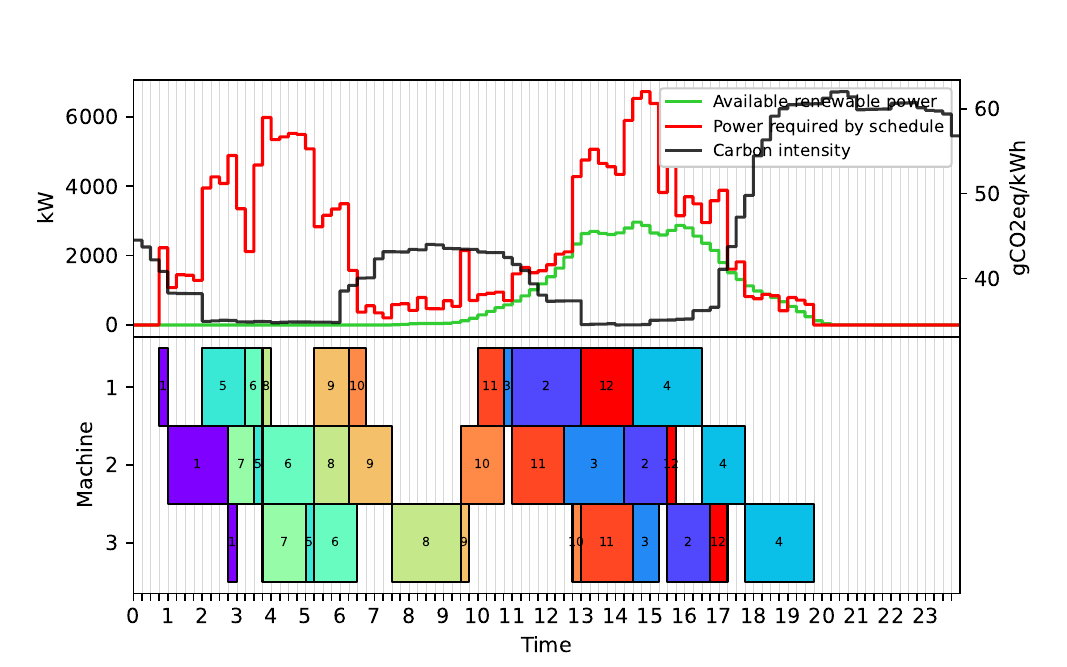}
		\caption{Carbon-aware scheduling: $4.85$ tCO\textsubscript{2}}
		\label{fig:gantt_CAS}
	\end{subfigure}
	\caption{Gantt chart, power requirement profiles, and resulting GHG emissions for the same instance optimized under two different objectives}
	\label{fig:makespan_vs_CAS}
\end{figure}

Figure \ref{fig:makespan_vs_CAS} illustrates two optimal schedules for instance \textit{CAS-PFSP-M3T1\_benchmark}: one optimized for makespan and the other for carbon emissions.
This instance consists of 12 jobs processed on 3 machines and the complete data are available in the public repository referenced in Section \ref{subsection:instance_datasets}.
While both achieve the same production goals, carbon-aware scheduling strategically delays certain operations to align energy-intensive processing with periods of low grid carbon intensity, thereby reducing overall emissions.

\subsection{Related work}
\label{subsection:related_work}

Algorithm configuration in EAs has been widely studied over the last three decades \cite{Eiben1999}.
Classically, a distinction is made between offline (static) algorithm configuration, where parameters are optimized prior to execution, and online (dynamic) configuration, where parameters are adjusted during the search.

Early work on online algorithm configuration primarily investigated parameter control strategies based on predefined mechanisms rather than learning.
In this line of work, parameter control methods are traditionally classified into three categories: deterministic, self-adaptive, and adaptive \cite{eiben2015introduction,Karafotias2015}.
Deterministic control methods modify parameters according to a predefined schedule, typically as a function of time or iteration count (e.g., \citet{Sun2020} and \citet{Flores-Torres2023}).
Self-adaptive control encodes parameters within the individuals' chromosomes, allowing parameter values to evolve via the same genetic operators that act on solutions (e.g., \citet{Back2000} and \citet{Gomez2004}).
Adaptive control, by contrast, relies on an external feedback mechanism that updates parameter values based on indicators of search progress such as population diversity, fitness improvement, or convergence rate (e.g., \citet{Thierens2005} and \citet{Aleti2013}).
For an extensive overview of work on parameter control in EAs, we refer the reader to \citet{Karafotias2015} and \citet{Aleti2016}.

Importantly, parameters in EAs can be either quantitative (e.g., mutation rate, population size) or qualitative, such as the selection among different variation operators \cite{Eiben2011}.
Within this broader perspective, parameter control can be viewed as a specific form of algorithm configuration.

More recently, Dynamic Algorithm Configuration \cite{Biedenkapp2020} has been proposed as a generalization of classical parameter control, modeling parameter adjustment as a sequential decision-making problem.
In DAC, parameter control is formulated as a Markov Decision Process (MDP) defined by a 4-tuple ($\mathcal{S}$, $\mathcal{A}$, $\mathcal{T}$, $\mathcal{R}$).
Here, $\mathcal{S}$ denotes the state space describing the current state of the algorithm (e.g., population statistics, search progress indicators), $\mathcal{A}$ the action space specifying possible parameter settings, $\mathcal{T}$ the transition distribution governing how the algorithm state changes after an action, and $\mathcal{R}$ a reward function quantifying the utility or progress of the search.
This formulation naturally enables the use of Reinforcement Learning \cite{SuttonBarto1998}, where an agent learns a parameter control policy through interaction with the running algorithm.

Several studies have applied RL-based DAC to EAs, with early work primarily focused on genetic algorithms (GAs).
Initial approaches mainly relied on tabular methods such as Q-learning \cite{Clifton2020} or SARSA \cite{Chen2005} to estimate an optimal parameter control strategy, which requires discretizing both the state and action spaces.
For instance, \citet{Chen2020} and \citet{Quevedo2021} dynamically adjusted two quantitative GA parameters, crossover rate and mutation probability, using Q-learning.
Although both studies adopted a similar control scheme, their applications differ: \citet{Chen2020} addressed the flexible job-shop scheduling problem (FJSP) \cite{DAUZEREPERES2024409} with makespan minimization, whereas \cite{Quevedo2021} considered the capacitated vehicle routing problem (CVRP) \cite{archetti2025}.
Extending this line of research, \citet{li2025} added a third, qualitative control parameter governing the selection among four local search operators, and applied their method to the charge scheduling problem for battery swap stations.
However, these approaches rely on tabular RL methods that require explicit enumeration of states and actions, which becomes impractical for EAs due to the high-dimensional and often continuous nature of the state and action spaces \cite{SuttonBarto1998}.

Consequently, more recent work has shifted towards RL with function approximation, in particular deep reinforcement learning (DRL) \cite{François-Lavet2018}, for learning parameter control policies.
By using multi-layer (i.e., deep) neural networks to approximate value functions, DRL alleviates the need for explicit state discretization and enables continuous, high-dimensional spaces \cite{Schuchardt2019}.
Several studies have therefore applied DRL to parameter control in differential evolution (DE), an EA proposed by \citet{storn1997} for continuous optimization.
Among these, \citet{Sharma2019} employed a Double Deep Q-Network agent (DDQN) \cite{vanhasselt2016} to select among four possible mutation strategies, while \citet{tan2021} used a DQN agent to choose between three alternatives.
\citet{sun2021} instead adopted a policy gradient approach \cite{SuttonBarto1998} to control two quantitative DE parameters: the scale factor and crossover rate.
Focusing on multi-objective optimization, \citet{Reijnen2022} combined these two parameters into a finite action set and learned an adaptive policy using Proximal Policy Optimization (PPO) \cite{Schulman2017}.
Continuing along this line, \citet{Tian2023} applied DQN to select among four possible operators for MOEA/D \cite{Zhang2007}, a widely used EA for multi-objective optimization.

Although the aforementioned studies provide evidence that learning-based DAC can improve algorithmic performance, they all evaluate their methods exclusively on CEC benchmark functions \cite{Garden2014}, whose representativeness for real-world combinatorial optimization problems remains debatable \cite{Piotrowski2023}.
To the best of our knowledge, only a limited number of studies have applied DRL-based DAC to real-world combinatorial optimization problems, among which \citet{Reijnen2023} represents one of the few examples.
However, their model was validated on only seven publicly available FJSP instances and reported no statistically significant improvement over a static GA with tuned parameters.
The authors concluded that, in this setting, the main benefit of DRL-based DAC lies in replacing parameter tuning with agent training, and emphasized the need for broader empirical validation.

As a result, it remains unclear whether the additional complexity and training cost of DRL-based DAC are justified for real-world optimization problems, especially when taking into account both generalization to unseen instance types and the computational cost of training.

\subsection{Scope of the paper}
\label{subsection:scope_of_the_paper}

Given these considerations, the objective of this paper is to critically assess the practical value of DRL-based DAC in a real-world optimization context.
Specifically, we propose a DRL-based DAC framework for CAS and design an experimental setup that explicitly evaluates when learning provides measurable benefits over tuning and when it does not.

In this context, our focus on DRL is motivated by its suitability for studying generalization.
DRL enables the learning of parameter control policies in high-dimensional and continuous spaces and can be trained in an instance-independent manner, making it well suited for investigating whether such control strategies generalize beyond the training distribution.
In contrast, classical parameter control approaches typically rely on predefined rules or problem-specific feedback mechanisms and often require manual design choices.
As a result, they are not designed with cross-instance generalization as a primary objective.
By focusing on DRL, we therefore aim to assess whether learning-based approaches can provide transferable control policies that remain effective beyond the training distribution.

To this end, we examine two complementary aspects.
First, we directly compare the proposed DRL-based approach with an otherwise identical static algorithm configuration, thereby isolating the contribution of learning from that of the underlying optimization algorithm.
Second, we assess performance across instance sets of increasing combinatorial complexity, allowing us to analyze how problem difficulty influences the usefulness of learned parameter control.

Building upon \citet{Reijnen2023}, the scope of this study is three-fold.
First, we evaluate the proposed framework on a substantially broader set of problem instances spanning different levels of combinatorial complexity.
Second, we explicitly investigate the generalizability of learned policies to previously unseen instances.
Third, we analyze under which conditions the additional modeling and computational effort associated with learning is justified by performance gains.

The remainder of the paper is structured as follows.
In Section \ref{section:method}, the proposed method is introduced.
Then, in Section \ref{section:computational_experiments}, the computational experiments and results are discussed.
Finally, conclusions and directions for future work are provided in Section \ref{section:conclusion}.

\section{Method}
\label{section:method}

The proposed MA-DRL-CAS-PFSP method builds upon the previously mentioned MA-CAS-PFSP introduced by \citet{Mencaroni2025}.
For ease of reading, we report in the section below the main characteristics of this algorithm.

\subsection{MA-CAS-PFSP Algorithm}
\label{subsection:MA-CAS-PFSP_algorithm}

The algorithm is a memetic algorithm (MA) that combines a GA with a local search operator to solve the carbon-aware permutation flow-shop scheduling problem.
As typical in GAs, the search is performed on a population of candidate solutions, which are also called individuals in EA literature.
Each individual encodes a complete production schedule using a dual random-key representation.
Specifically, solutions consist of (i) a job sequence determining the processing order on all machines and (ii) machine-specific idle-time allocations that distribute available slack time within the planning horizon.
This representation allows the algorithm to jointly optimize job sequencing and the temporal placement of idle time, so to align energy-intensive operations with periods of lower carbon intensity of electricity.

During the search process, a number of evolutionary operations are performed on the individuals.
First, offspring are generated through a controlled swap crossover operator, which exchanges subsets of genes between two parents.
A nonuniform mutation operator is then applied to the offspring that perturbs random keys with adjustable frequency and severity.
To intensify the search, each offspring solution is further refined using a local search procedure based on pairwise adjacent job swaps.
Selection is elitist: at each generation, the best-performing individuals are retained for the next iteration.
The algorithm terminates after a predetermined number of generations.

Several parameters govern the behavior of the algorithm.
These parameters can be organized into three families, reflecting their functional role in the search process.

The first family comprises generic EA parameters, which define the overall computational budget.
These include the population size $\rho$ and the maximum number of generations $\gamma$.
In this study, these parameters are kept fixed.

The second family consists of crossover-related parameters, which control how information is recombined between parent solutions.
A fraction $\xi$ of the new population is produced via crossover.
Crossover itself operates independently on the two components of the representation (i.e., job sequence and idle-time allocation).
The parameters $\chi_j$ and $\chi_p$ regulate the severity of crossover by specifying, for each gene, the probability that it is swapped between parents in the job-sequence and pause arrays, respectively.

The third family contains mutation-related parameters, which govern the stochastic perturbation of solutions.
Similarly to crossover, mutation is applied independently to the job sequence and pause components.
The parameters $\pi_j$ and $\pi_p$ control the frequency of mutation by expressing the probability that a given gene is mutated.
The severity of mutation is determined by $\sigma_j$ and $\sigma_p$ by setting the standard deviation of the normal perturbation applied to the corresponding random keys.

Taken together, the crossover and mutation parameters allow the algorithm to regulate both the frequency and the intensity of variation.
In the static configuration, these parameters are tuned offline and kept fixed throughout the search.
Details on the tuning procedure and the resulting parameter values are provided in Section \ref{subsection:static_parameter_tuning}.

In the remainder of this paper, this algorithm serves as the reference point against which we evaluate the proposed DRL-based DAC.
For a detailed description, we refer the reader to \citet{Mencaroni2025}.

\subsection{Markov Decision Process Formulation}
\label{subsection:MDP_formulation}

In order to apply DRL, we need to formulate the DAC problem as an MDP.
This involves defining the four components that fully describe the process at each iteration of the algorithm.
In what follows, each component of the MDP formulation is detailed, including the state representation used as input to the agent and the reward mechanism that drives the learning process.

\subsubsection{States}
\label{subsubsection:states}
The state space provides the agent with information about the current status and progress of the search.
In our formulation, the state is represented by a five-dimensional vector of normalized metrics, which serves as input to the DRL agent: current best fitness, mean fitness, coefficient of variation, remaining budget, and stagnation count.
This representation is designed to include all and only the information necessary to characterize the search dynamics in an instance-independent way, enabling the agent to learn policies that can generalize across diverse instance types.

The first two components correspond to the best and average fitness values (i.e., lowest objective value) in the current population and provide a direct indication of the search progress.
Both metrics are normalized with respect to the best fitness value of the initial generation in order to provide instance-independent observations and enable learning across heterogeneous instances.

The coefficient of variation, defined as the ratio of the standard deviation to the mean fitness value of the current population, serves as a proxy for the diversity of the population.
The remaining budget captures the fraction of iterations left before termination, while the stagnation count represents the normalized number of consecutive generations without improvement of the best-found solution.

All state components are normalized to the $[0,\,1]$ range to ensure instance independence, which supports the generalizability of learned policies across problem instances.

\subsubsection{Actions}
\label{subsubsections:actions}
The action space is defined as a seven-dimensional continuous vector, with each component taking values in the interval $[-1,\,1]$.
Each action directly controls one of the configurable parameters of the MA-CAS-PFSP algorithm introduced in Section \ref{subsection:MA-CAS-PFSP_algorithm}.
An action, consisting of seven parameter values, is selected once per generation.
These values are used for one iteration of the algorithm, until a new population is generated.

To facilitate learning, actions are linearly rescaled to predefined parameter-specific ranges, as reported in Table \ref{tab:rescaling_actions}.
Although not strictly required, bounding the action space has been shown to improve learning efficiency by preventing the exploration of implausible configurations \cite{Huang2022}.
The selected rescaling ranges correspond to commonly used values in EA literature \cite{CoelloCoello2007}, and were chosen sufficiently broad to avoid constraining the policy search space.

\begin{table}[h]
   	\fontsize{9pt}{14pt}\selectfont
   	\centering
       \caption{Rescaling limits of algorithm parameters for action space}
   	\begin{tabular}{c|c|c}
   		Parameter & Minimum value & Maximum value\\
   		\hline
   		$\xi$ & 0.5 & 0.9 \\
   		$\chi_j$ & 0.1 & 0.5 \\
   		$\chi_p$ & 0.05 & 0.5 \\
   		$\pi_j$ & 0.01 & 0.2 \\
   		$\pi_p$ & 0.01 & 0.11 \\
   		$\sigma_j$ & 0.008 & 0.2 \\
   		$\sigma_p$ & 0.15 & 0.25 \\
   	\end{tabular}
   	\label{tab:rescaling_actions}
\end{table}

\subsubsection{Transitions}
\label{subsubsection:transitions}

The state transition function describes how the environment evolves after the agent applies an action.
In principle, this function maps the current state and the selected algorithm parameters to the next state of the optimization process.
However, in the present setting, the transition dynamics arise from the complex and stochastic interaction between evolutionary operators, selection mechanisms, and the underlying scheduling problem.
As a result, the transition function cannot be expressed in closed form nor derived from the algorithmic implementation \cite{Biedenkapp2020}.

Therefore, the transition dynamics are not modeled explicitly within the MDP.
Instead, the agent learns an implicit model of the environment through repeated interaction with the MA-CAS-PFSP algorithm, observing the resulting state transitions induced by its parameter control actions.

\subsubsection{Rewards}
\label{subsubsection:rewards}

To enable the agent to learn an effective dynamic parameter control policy, a reward signal is given at each interaction step.
Following prior work on RL for DAC \cite{Biedenkapp2020}, we adopt a dense reward scheme, whereby the agent receives feedback after every generation of the underlying EA.

The reward is designed to reflect the improvements in solution quality achieved between consecutive generations.
However, EAs typically exhibit rapid progress during the early stages of the search, while later stages require increased effort to yield even marginal improvements \cite{eiben2015introduction}.
To account for this behavior, we adopt a reward function formulation inspired by \citet{Reijnen2025}, which places increasing emphasis on improvements achieved in later, more challenging phases of the search.

At each iteration $t$, the reward is computed based on the change in the best fitness within the population.
If the best fitness in the current iteration $f_{current}$ improves upon the best fitness observed so far $f_{previous}$, the reward is defined as the difference between the squared normalized improvements.
If no improvement is observed, the reward is set to zero.
Formally, the reward at iteration $t$ is given by:
\begin{equation}
    r_t = \begin{cases}
        {\delta^2_{current} - \delta^2_{previous}} & \text{if} \enspace f_{current} < f_{previous}\\
        {0} & \text{otherwise}
    \end{cases} \enspace .
\end{equation}
The normalized improvements $\delta_{current}$ and $\delta_{previous}$ are defined as:
\begin{equation}
    \delta_{current} = \left( \frac{f_{initial} - f_{current}}{f_{initial} - f_{ideal}} \right) \cdot 100 \enspace ,
\end{equation}
and
\begin{equation}
    \delta_{previous} = \left( \frac{f_{initial} - f_{previous}}{f_{initial} - f_{ideal}} \right) \cdot 100 \enspace .
\end{equation}

Here, $f_{initial}$ denotes the best fitness found in the first generation, while $f_{ideal}$ corresponds to the best solution obtained by an extended offline run of the static MA-CAS-PFSP algorithm, in which the computational budget is doubled.
These reference values normalize the reward signal, making it proportional to the effective progress achieved relative to both the starting point of the search and a high-quality, yet realistically attainable, target solution obtained from the same algorithmic framework.

This design provides a normalized reward signal that is comparable across instances, which is essential for learning policies that generalize beyond the training set.
While the use of a reference solution obtained from extended runs of the static algorithm may introduce a bias toward regions of the search space explored by this algorithm, it allows the agent to evaluate improvements relative to a strong and achievable baseline, without requiring access to optimal solutions, which are computationally infeasible to obtain at the considered scale.

Importantly, the reward function is defined in terms of relative improvement between consecutive iterations, rather than similarity to the reference solution itself.
As a result, the agent is rewarded for making progress in the search process, irrespective of the specific trajectory followed.
In particular, if the agent discovers solutions that outperform $f_{ideal}$, the normalized improvement exceeds 100, resulting in larger rewards.
This further encourages improvements beyond the reference solution and highlights that the reward function does not impose a strict upper bound on achievable performance.

\begin{figure}[h]
	\centering
	\includegraphics[width=.8\linewidth]{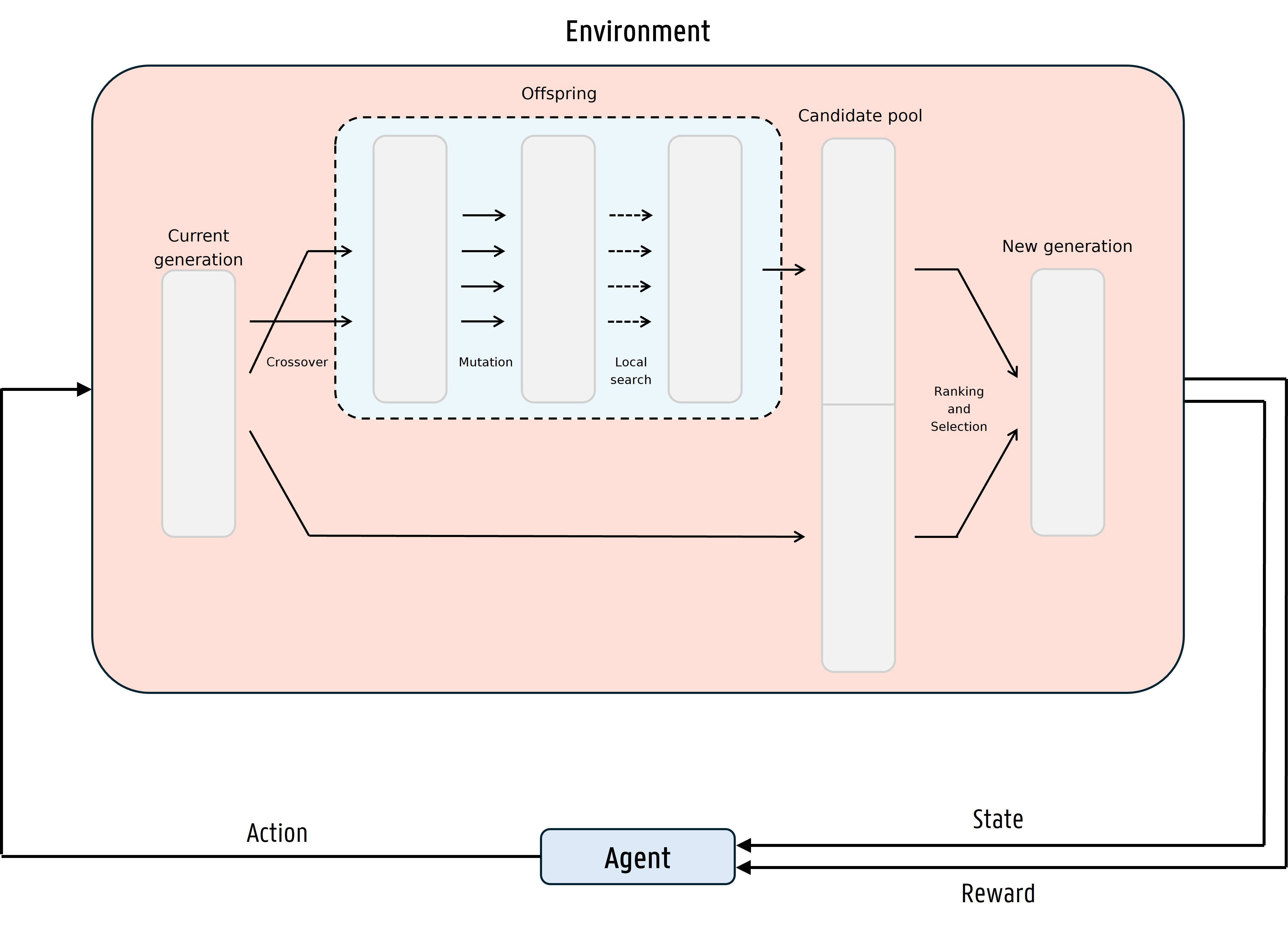}
	\caption{Overview of the MA-DRL-CAS-PFSP algorithm}
	\label{fig:overview_MA-DRL}
\end{figure}

A schematic overview of the proposed MA-DRL-CAS-PFSP algorithm is shown in Figure \ref{fig:overview_MA-DRL}.

\section{Computational experiments}
\label{section:computational_experiments}

To assess the generalization capabilities of the proposed method, we conducted a series of computational experiments.
In this section, we describe the experimental setup used and present the corresponding results.

\subsection{Instance datasets}
\label{subsection:instance_datasets}
We consider ten instance datasets covering a range of combinatorial complexity.
Each dataset consists of fifty randomly generated instances of similar characteristics and complexity.

To assess generalizability, we distinguish between \textit{known} and \textit{unknown} instance types.
The former comprises the four datasets introduced in \citet{Mencaroni2025}, which are used during training.

The latter consists of six additional datasets generated using the same instance-generation procedure, but exhibiting higher combinatorial complexity, primarily through an increased number of machines.
In addition to the increase in machine count, the number of jobs also varies across datasets, although to a lesser extent, as reflected in Table \ref{tab:overview_instance_datasets}.

These datasets are used exclusively for evaluation.
Since these instance types are not observed during training, we consider these instances as unknown, and evaluate in our experiments how well the algorithms perform on such unseen instances.
This separation enables a controlled analysis of the generalization capabilities of the learned parameter control policies.

An overview of all considered datasets along with their key characteristics is provided in Table \ref{tab:overview_instance_datasets}.
The full instance files are available at \url{https://github.com/ugent-isye/MA-DRL-CAS-PFSP}. 

\begin{table}[h]
	\fontsize{9pt}{14pt}\selectfont
	\centering
    \caption{Overview of the instance datasets and their key features}
	\begin{tabular}{c|c|c|c|c|c|c}
		Name & \# instances & $M$ & $T$ & \# jobs & \# operations & Type \\
		\hline
		CAS-PFSP-M1T1 & $50$ & $1$ & $1$ & 6 - 15 & 6 - 15 & Known \\
		CAS-PFSP-M1T3 & $50$ & $1$ & $3$ & 25 - 40 & 25 - 40 & Known \\
		CAS-PFSP-M3T1 & $50$ & $3$ & $1$ & 8 - 18 & 24 - 54 & Known \\
		CAS-PFSP-M3T3 & $50$ & $3$ & $3$ & 34 - 61 & 102 - 183 & Known \\
        CAS-PFSP-M5T1 & $50$ & $5$ & $1$ & 5 - 13 & 25 - 65 & Unknown \\
        CAS-PFSP-M5T3 & $50$ & $5$ & $3$ & 29 - 51 & 145 - 255 & Unknown \\
        CAS-PFSP-M10T1 & $50$ & $10$ & $1$ & 13 - 20 & 130 - 200 & Unknown \\
        CAS-PFSP-M10T3 & $50$ & $10$ & $3$ & 65 - 83 & 650 - 830 & Unknown \\
        CAS-PFSP-M15T1 & $50$ & $15$ & $1$ & 23 - 31 & 345 - 465 & Unknown \\
        CAS-PFSP-M15T3 & $50$ & $15$ & $3$ & 106 - 139 & 1590 - 2085 & Unknown \\
	\end{tabular}
	\label{tab:overview_instance_datasets}
\end{table}

\subsection{Training procedure}
\label{subsection:training_procedure}

In line with \citet{Reijnen2023}, the parameter control policy was trained using Proximal Policy Optimization (PPO), a state-of-the-art policy gradient algorithm.
Training was conducted on a fixed set of representative training instances from the known datasets introduced in \citet{Mencaroni2025}.
Specifically, three instances were selected from each of the four datasets, resulting in a total of 12 training instances.
These instances were kept separate from the test instances used in the computational experiments.

The selected training instances are representative of their respective datasets, which are generated using a common procedure and therefore share similar structural characteristics.
This allows the training process to capture generalizable patterns at the dataset level, rather than memorizing instance-specific behavior.

Each training episode was initialized by selecting a training instance uniformly at random and executing 100 iterations of the underlying MA-CAS-PFSP algorithm, using a population size of 250 individuals.
During an episode, the agent interacted online with the algorithm by observing the current state, adapting its parameter configuration, and receiving rewards proportionally to the effectiveness of its actions.

The total duration of the training process was determined empirically based on the observed convergence behavior of the learning curve.
Training was terminated once the marginal improvements in the mean episode reward became negligible, resulting in a total of four million training steps, where each step corresponds to one iteration of the algorithm applied to a single instance.
On the hardware used in this study, this corresponded to an average processing speed of approximately 20 steps per second, resulting in a total training time of about 55 hours (approximately 2.3 days).
We note that the exact wall-clock time depends on hardware and implementation details.

\begin{figure}[H]
	\centering
	\includegraphics[width=.6\linewidth, trim={0cm 0cm 0cm 0cm}, clip]{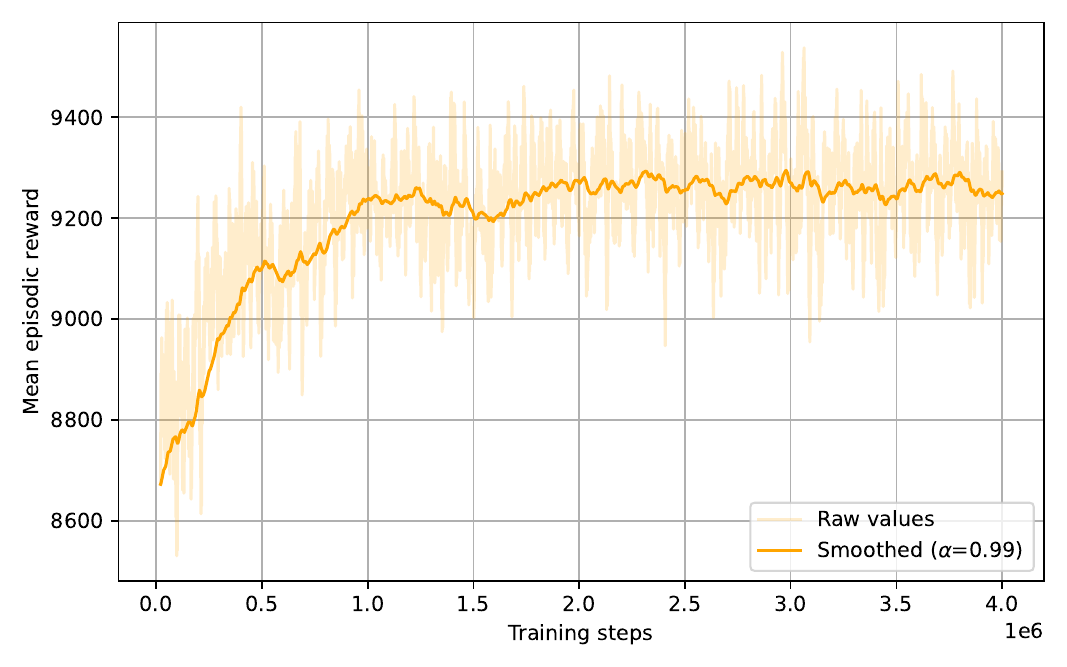}
	\caption{Mean reward per episode during training}
	\label{fig:training_plot}
\end{figure}

Figure \ref{fig:training_plot} reports the evolution of the mean episode reward over the course of the training, showing a clear convergence trend.
The PPO agent was implemented using the Stable-Baselines3 library with default hyperparameters\footnote{\url{https://stable-baselines3.readthedocs.io/en/master/modules/ppo.html\#parameters}}, including a learning rate of $3 \cdot 10^{-4}$, a discount factor $\gamma = 0.99$, a clipping range of 0.2, and a batch size of 64. The policy network consists of a multilayer perceptron with two hidden layers of 64 neurons each and ReLU activation functions.

\subsection{Static parameter tuning}
\label{subsection:static_parameter_tuning}

For static parameter tuning, we employed an automated configuration approach based on Optuna \cite{Optuna2019}, and additionally validated the robustness of the resulting configuration using irace \cite{LopezIbanez2016}, both under the same computational budget and the same set of training instances.

To enable a fair comparison with the proposed DRL-based approach, the tuning setup was designed to mirror the DRL training setup as closely as possible.
Specifically, tuning was performed on the same set of 12 training instances used during agent training, and the total computational budget allocated to tuning was matched to the training budget of four million algorithm iterations.

In this setting, one algorithm iteration corresponds to executing one generation of the memetic algorithm on a single instance.
Each candidate configuration is therefore evaluated by running the algorithm on all 12 training instances for 100 generations with a population size of 250 individuals.
Since the population size is fixed, the computational effort is expressed in terms of generations, resulting in $12 \cdot 100 = 1\,200$ algorithm iterations per configuration evaluation.
Given the total budget of four million algorithm iterations, this corresponds to approximately $3\,333$ configuration evaluations.
During tuning, the effectiveness of each configuration is assessed by calculating the average objective value across the 12 training instances.

Optuna explores the configuration space using a Tree-structured Parzen Estimator strategy, while irace employs an iterative racing procedure to discard poorly performing configurations and focus the search on promising regions.
For both methods, the configuration yielding the lowest average objective value is selected.

Under this setup, the total wall-clock time for tuning was approximately 3 days on the same hardware.
We note that runtime is influenced by implementation details and the degree of parallelization, and is therefore reported for reference only.

The configurations obtained via Optuna and irace were found to yield highly comparable performance across all datasets, with only negligible differences and no consistent advantage for either method.
This indicates that the tuned configuration is robust with respect to the choice of tuning procedure and confirms that the results of this study are not sensitive to the specific configuration method used.

An overview of the algorithm parameters, together with their default and the configurations obtained via Optuna (used in the main experiments) and irace (used for validation), is provided in Table \ref{tab:algorithm_parameters_overview}.

\begin{table}[h]
   	\fontsize{9pt}{14pt}\selectfont
   	\centering
       \caption{Overview of algorithm parameters with their default and tuned values}
   	\begin{tabular}{c|l|c|c|c}
   		Parameter & Definition & Default & Tuned (Optuna) & Tuned (irace)\\
   		\hline
   		$\rho$ & Population size & 250 & 250 & 250 \\
   		$\gamma$ & Maximum number of generations & 100 & 100 & 100 \\
   		$\xi$ & Crossover rate & 0.85 & 0.885 & 0.780 \\
   		$\chi_j$ & Crossover probability for jobs & 0.5 & 0.418 & 0.490 \\
   		$\chi_p$ & Crossover probability for pauses & 0.5 & 0.133 & 0.054 \\
   		$\pi_j$ & Mutation probability for jobs & 0.05 & 0.017 & 0.020 \\
   		$\pi_p$ & Mutation probability for pauses & 0.05 & 0.014 & 0.015 \\
   		$\sigma_j$ & Mutation step size for jobs & 0.2 & 0.012 & 0.039 \\
   		$\sigma_p$ & Mutation step size for pauses & 0.2 & 0.233 & 0.248 \\
   	\end{tabular}
   	\label{tab:algorithm_parameters_overview}
\end{table}

\subsection{Experiments design}
\label{subsection:experiments_design}

We conducted ten independent sets of experiments, one for each dataset introduced in Section \ref{subsection:instance_datasets}.
For every dataset, its 50 problem instances were solved 10 times using each of the following three algorithm variants: (i) the static MA algorithm with default parameter setting, (ii) the static MA with tuned parameter setting, and (iii) the proposed MA-DRL algorithm.

Repeating each instance ten times accounts for the inherent stochasticity of evolutionary search.
For each dataset-method combination, we first compute the mean and the best objective value across the ten runs.
To compensate for variability in difficulty between instances, results are then aggregated at the dataset level.
Specifically, for each dataset-method combination, we report the average (over the 50 instances) of the mean objective values computed over the 10 runs, as well as the average of the best objective values obtained over the 10 runs.

To assess the performance of the algorithms, we report the relative percentage improvement $\% \Delta$ of the DRL-based approach with respect to the static baseline, defined as:
\begin{equation}
	\% \Delta = \frac{\left( \text{Obj - Obj\textsubscript{MA-DRL}} \right)}{\text{Obj}} \cdot 100 \enspace .
\end{equation}
Accordingly, a positive value of $\% \Delta$ indicates that the DRL-based parameter control approach outperforms the static algorithm.

All the experiments were conducted on a machine equipped with an 11th Gen Intel(R) Core(TM) i5-1145G7 processor (2.60 GHz), using Python version 3.10.0.rc2.

\subsection{Results}
\label{subsection:results}

We now report the outcomes of the computational experiments described in Section \ref{subsection:experiments_design}, distinguishing between known and unknown instance types.

\subsubsection{Known instance types}

The first set of experiments concerns the known instance types, i.e., datasets whose characteristics closely resemble the instances used during the training of the DRL agent and the static parameter tuning.
These experiments therefore evaluate the behavior of the algorithms in deployment scenarios that are well represented during learning and tuning.

\begin{figure}[h]
	\centering
		\begin{subfigure}{.49\columnwidth}
			\centering
			\includegraphics[clip, trim={0cm 0cm 0cm 0cm},width=\textwidth]{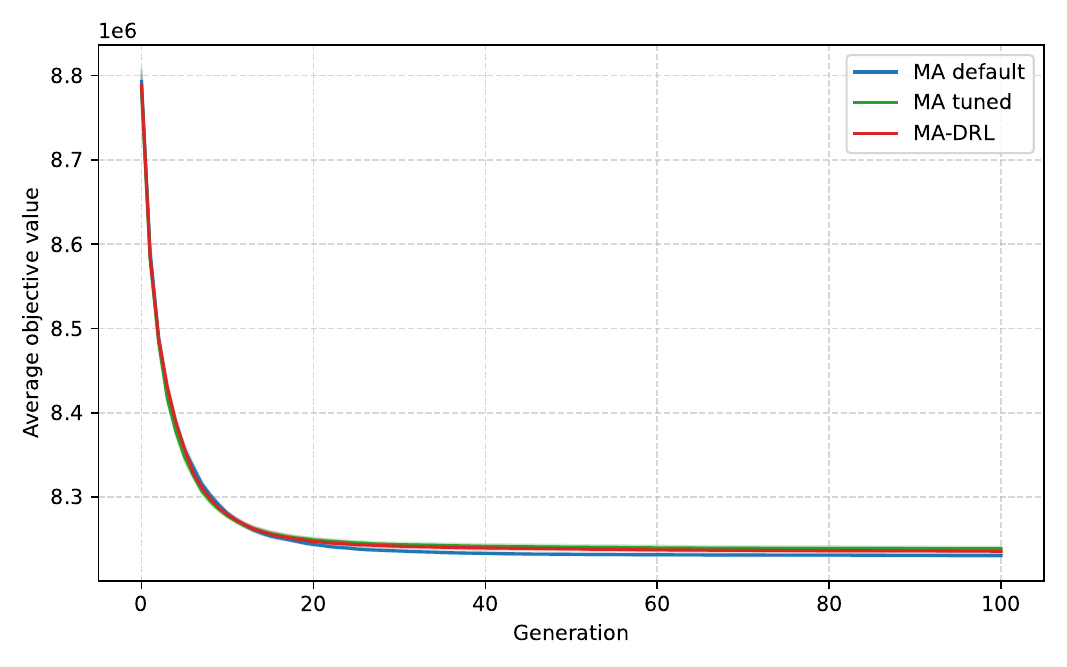}
			\caption{CAS-PFSP-M1T1}
			\label{fig:results_convergence_plot_M1T1}
		\end{subfigure}
		\begin{subfigure}{.49\columnwidth}
			\centering
			\includegraphics[clip, trim={0cm 0cm 0cm 0cm},width=\textwidth]{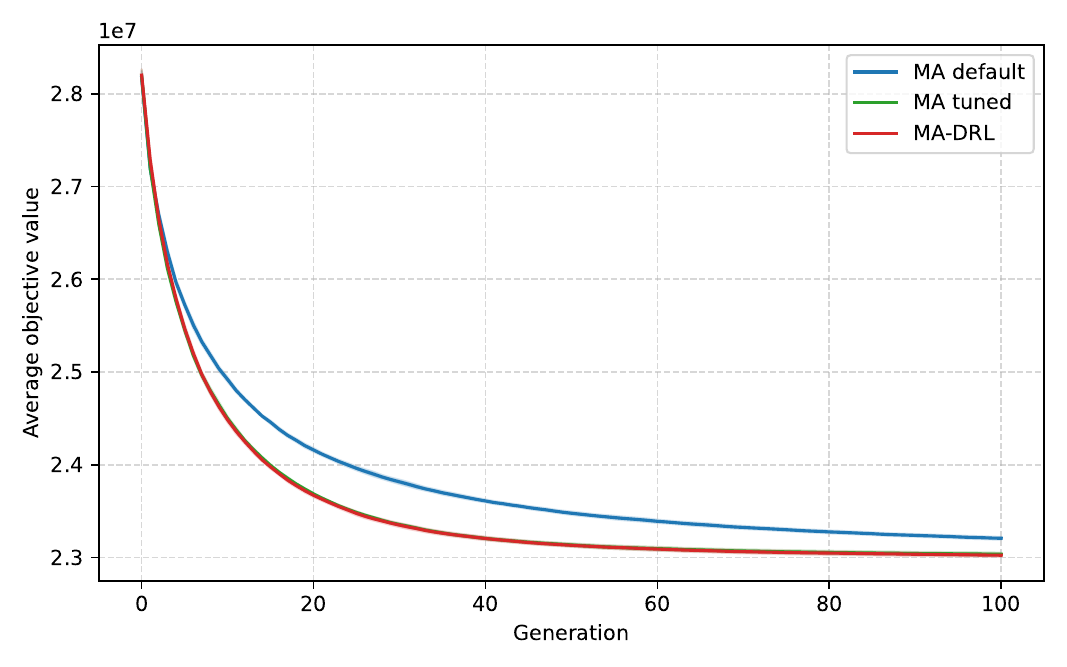}
			\caption{CAS-PFSP-M1T3}
			\label{fig:results_convergence_plot_M1T3}
		\end{subfigure}
		\begin{subfigure}{.49\columnwidth}
			\centering
			\includegraphics[clip, trim={0cm 0cm 0cm 0cm},width=\textwidth]{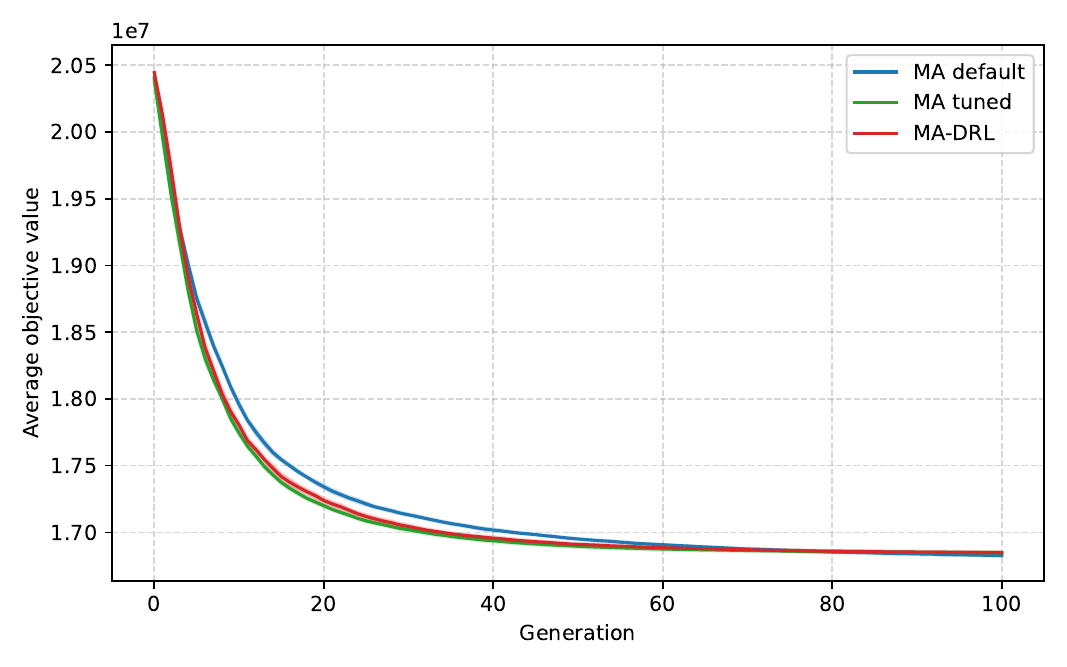}
			\caption{CAS-PFSP-M3T1}
			\label{fig:results_convergence_plot_M3T1benchmark}
		\end{subfigure} 
		\begin{subfigure}{.49\columnwidth}
			\centering
			\includegraphics[clip, trim={0cm 0cm 0cm 0cm},width=\textwidth]{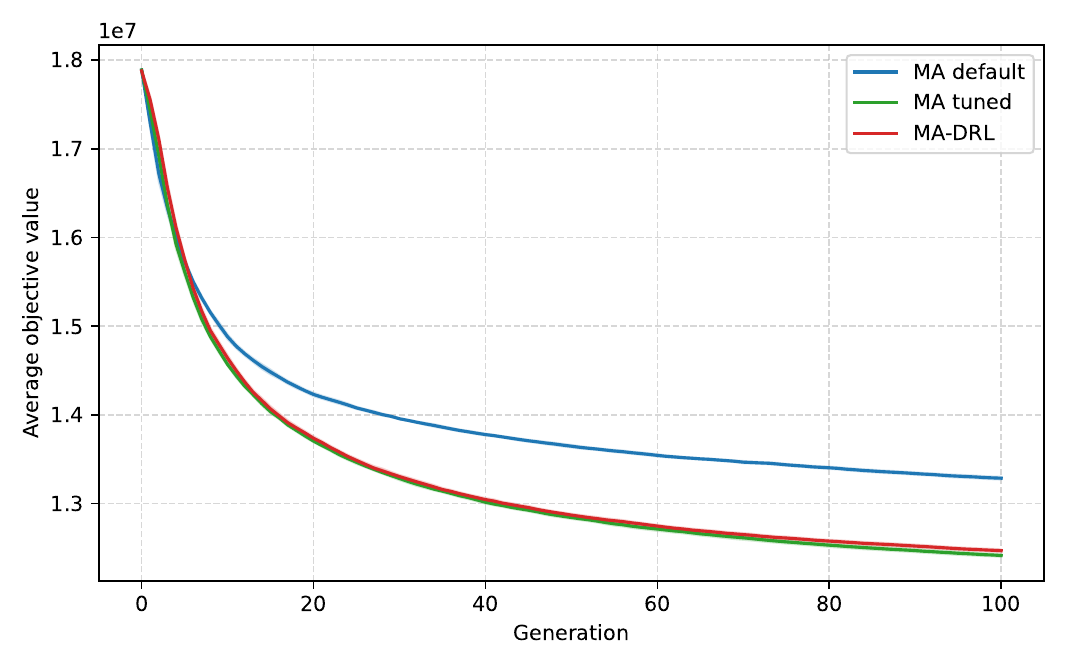}
			\caption{CAS-PFSP-M3T3}
			\label{fig:results_convergence_plot_M3T3}
		\end{subfigure}
	\caption{Convergence plots of average objective value per dataset, known instance types}
	\label{fig:results_convergence_plot_all_known}
\end{figure}

Table \ref{tab:overview_results_known} summarizes the numerical results for all four known datasets, whereas Figure \ref{fig:results_convergence_plot_all_known} visualizes the convergence behavior of the algorithms.
The table reports the mean, best, and standard deviation of the objective values across all instances, as well as the relative percentage improvements with respect to the DRL-based method.
For each dataset, pairwise comparisons were conducted between the best-performing method and the remaining approaches.
Best results are highlighted in bold, and underlined entries denote statistical significance according to the Wilcoxon rank sum test $(p < 0.05)$, with Holm correction applied to account for multiple comparisons against the best-performing method per dataset.

\begin{table}[h]
	\fontsize{9pt}{14pt}\selectfont
	\centering
    \caption{Overview of the results, known instance types}
	\begin{tabular}{c|c|c|c|c|c|c|c|c}
        \hline
        & \multicolumn{4}{c|}{CAS-PFSP-M1T1} & \multicolumn{4}{c}{CAS-PFSP-M1T3} \\
        \hline
        Method & mean & best & std & $\Delta\%$ & mean & best & std & $\% \Delta$ \\
        \hline
        MA default & $\mathbf{\underline{2.06 \times 10^6}}$ & $\mathbf{2.05 \times 10^6}$ & $2.78 \times 10^3$ & -0.05 & $5.80 \times 10^6$ & $5.73 \times 10^6$ & $4.49 \times 10^4$ & 0.70 \\
        MA tuned & $2.06 \times 10^6$ & $2.06 \times 10^6$ & $3.34 \times 10^3$ & 0.04 & $\mathbf{5.76 \times 10^6}$ & $\mathbf{5.69 \times 10^6}$ & $4.22 \times 10^4$ & -0.04 \\
        MA-DRL & $2.06 \times 10^6$ & $2.06 \times 10^6$  & $2.62 \times 10^3$ & 0 & $5.76 \times 10^6$ & $5.69 \times 10^6$ & $4.26 \times 10^4$ & 0 \\
        \hline
        & \multicolumn{4}{c|}{CAS-PFSP-M3T1} & \multicolumn{4}{c}{CAS-PFSP-M3T3} \\
        \hline
        Method & mean & best & std & $\Delta\%$ & mean & best & std & $\% \Delta$ \\
        \hline
        MA default & $\mathbf{\underline{4.20 \times 10^6}}$ & $\mathbf{4.18 \times 10^6}$ & $1.75 \times 10^4$ & -0.23 & $1.33 \times 10^7$ & $1.30 \times 10^7$ & $1.40 \times 10^5$ & 6.53 \\
        MA tuned & $4.22 \times 10^6$ & $4.18 \times 10^6$ & $2.55 \times 10^4$ & 0.05 & $\mathbf{\underline{1.24 \times 10^7}}$ & $\mathbf{1.22 \times 10^7}$ & $1.36 \times 10^5$ & -0.28 \\
        MA-DRL & $4.21 \times 10^6$ & $4.18 \times 10^6$  & $2.45 \times 10^4$ & 0 & $1.25 \times 10^7$ & $1.22 \times 10^7$ & $1.54 \times 10^5$ & 0 \\
        \hline
    \end{tabular}
    \label{tab:overview_results_known}
\end{table}

Across these datasets, the three methods exhibit remarkably similar performance.
In most cases, the numerically best-performing method is also statistically significantly better, yet the magnitude of these improvements is extremely small.
Although the significance tests indicate differences beyond random fluctuation, their practical impact is negligible.

In none of the four datasets does the DRL-based approach provide a measurable improvement over the tuned baseline.
Its performance is consistently close, occasionally marginally worse and occasionally marginally better, but never in a way that would suggest a practical advantage.
For example, in M1T1, MA-DRL improves over the tuned static method by only $0.04 \, \%$, while for M3T3 the tuned baseline performs $0.28 \, \%$ better than MA-DRL.

Furthermore, for the M1T1 and M3T1 datasets, tuning itself does not lead to a meaningful improvement over the default parameter setting.
Only in the more challenging M1T3 and M3T3 datasets, where the planning horizon is longer, does tuning provide a modest benefit over the default configuration, yet even there the tuned and the DRL-based method perform very similarly.

Overall, for the known instance types, we conclude that the DRL-based DAC does not offer an advantage over static tuning.
When the characteristics of new problem instances closely match those observed during training and tuning, static methods already achieve highly effective parameter configurations, and the DRL policy essentially replicates this performance without surpassing it.

\subsubsection{Unknown instance types}

\begin{figure}[h]
	\centering
		\begin{subfigure}{.49\columnwidth}
			\centering
			\includegraphics[clip, trim={0cm 0cm 0cm 0cm},width=\textwidth]{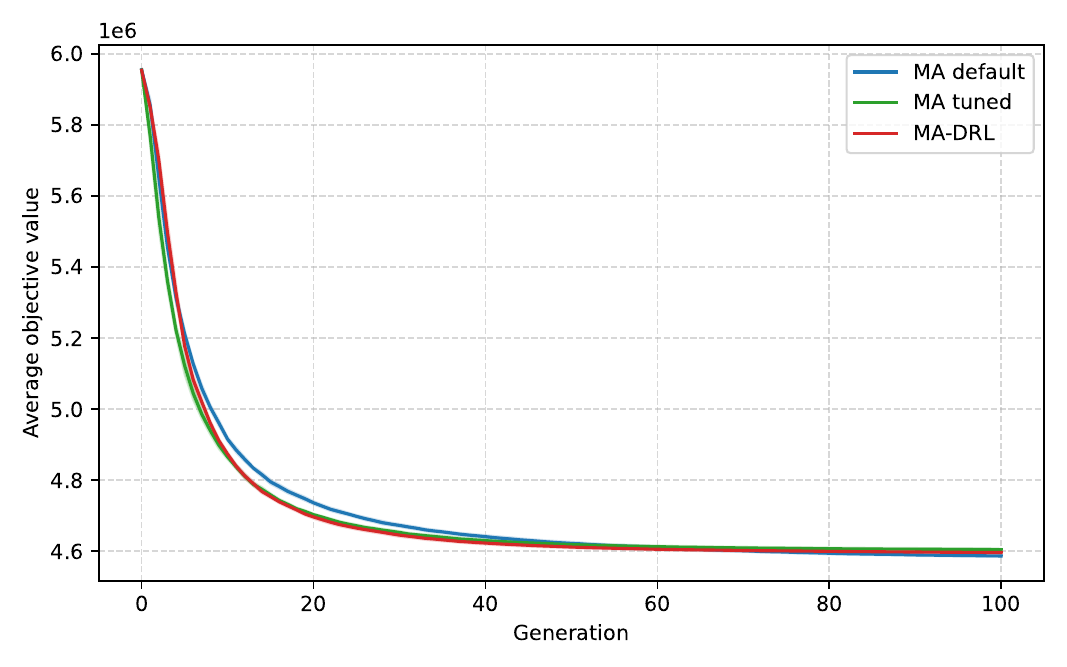}
			\caption{CAS-PFSP-M5T1}
			\label{fig:results_convergence_plot_M5T1}
		\end{subfigure}
		\begin{subfigure}{.49\columnwidth}
			\centering
			\includegraphics[clip, trim={0cm 0cm 0cm 0cm},width=\textwidth]{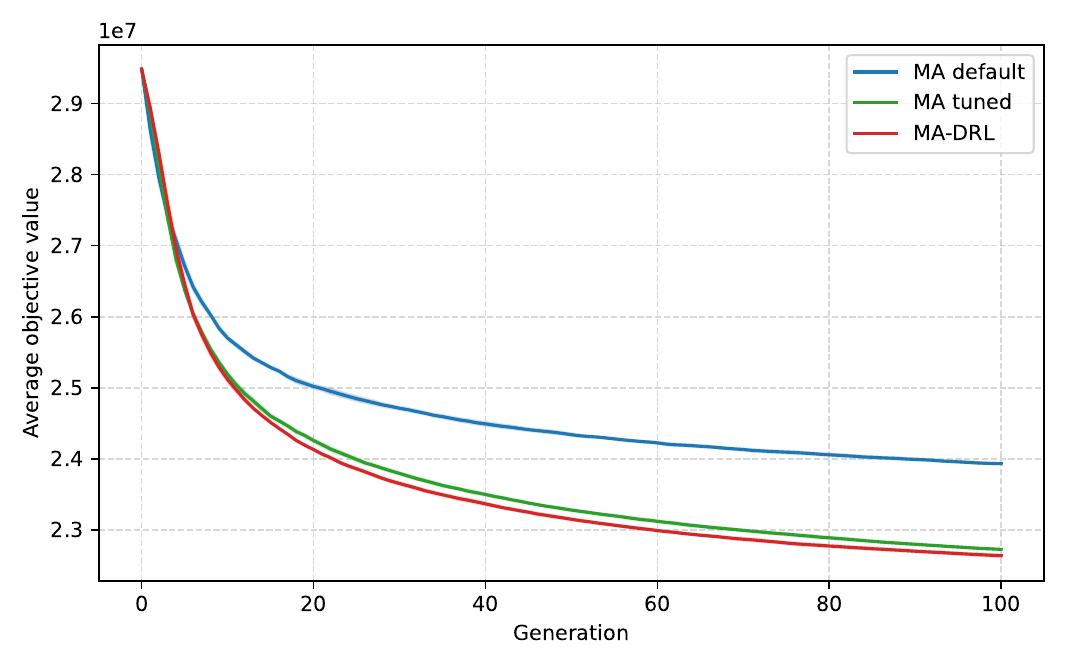}
			\caption{CAS-PFSP-M5T3}
			\label{fig:results_convergence_plot_M5T3}
		\end{subfigure}
		\begin{subfigure}{.49\columnwidth}
			\centering
			\includegraphics[clip, trim={0cm 0cm 0cm 0cm},width=\textwidth]{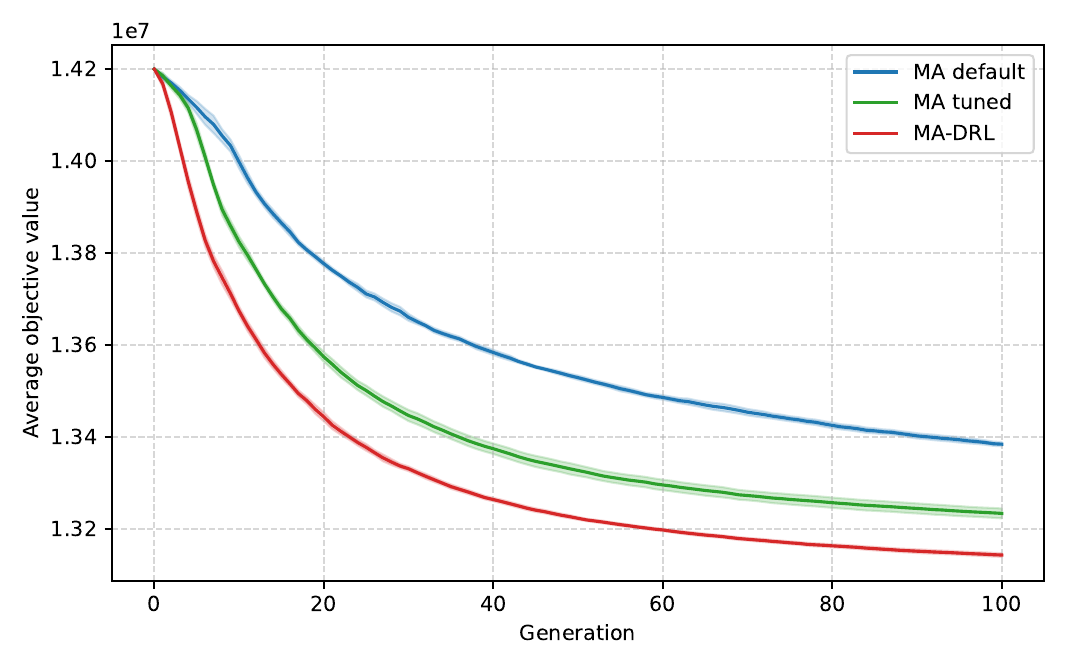}
			\caption{CAS-PFSP-M10T1}
			\label{fig:results_convergence_plot_M10T1}
		\end{subfigure}
		\begin{subfigure}{.49\columnwidth}
			\centering
			\includegraphics[clip, trim={0cm 0cm 0cm 0cm},width=\textwidth]{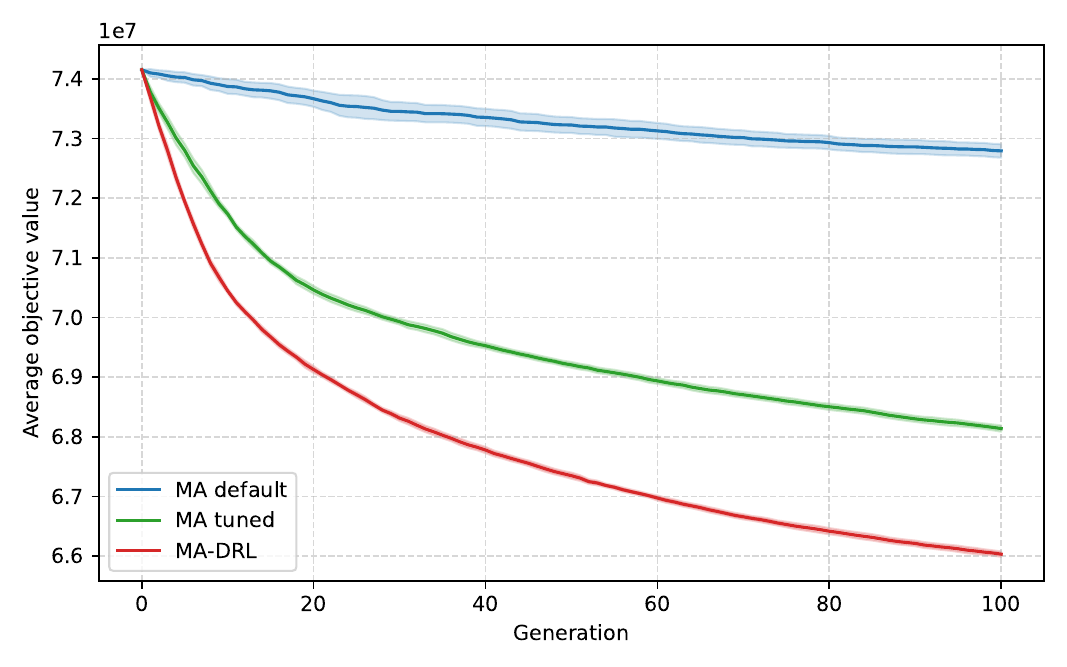}
			\caption{CAS-PFSP-M10T3}
			\label{fig:results_convergence_plot_M10T3}
		\end{subfigure}
		\begin{subfigure}{.49\columnwidth}
			\centering
			\includegraphics[clip, trim={0cm 0cm 0cm 0cm},width=\textwidth]{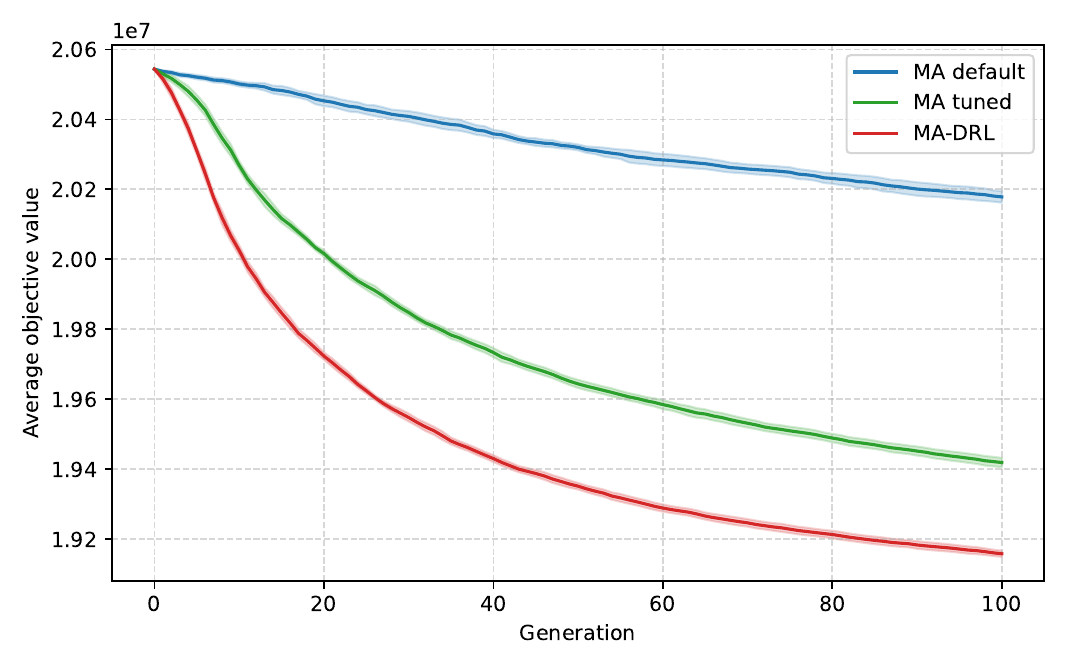}
			\caption{CAS-PFSP-M15T1}
			\label{fig:results_convergence_plot_M15T1}
		\end{subfigure}
		\begin{subfigure}{.49\columnwidth}
			\centering
			\includegraphics[clip, trim={0cm 0cm 0cm 0cm},width=\textwidth]{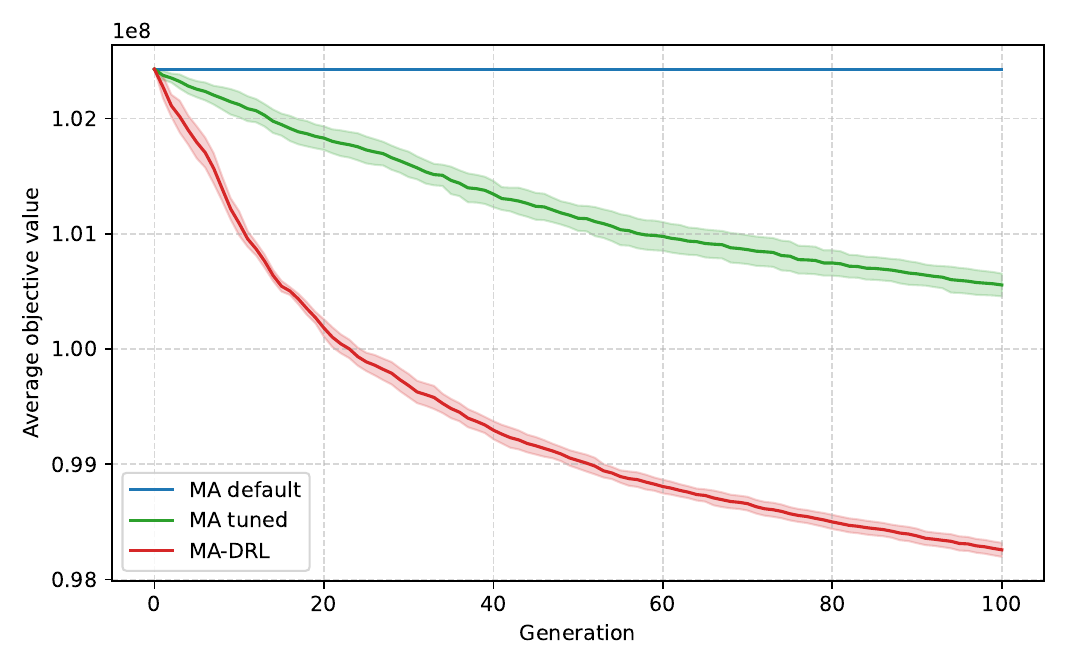}
			\caption{CAS-PFSP-M15T3}
			\label{fig:results_convergence_plot_M15T3}
		\end{subfigure}
	\caption{Convergence plots of average objective value per dataset, unknown instance types}
	\label{fig:results_convergence_plot_all_unknown}
\end{figure}

We now turn to the unknown instance types, which are six datasets generated with the same procedure as the known datasets but of higher combinatorial complexity, primarily due to a larger number of machines.
These instance types were never observed during training or tuning and thus provide a basis for evaluating the generalization capability of the DRL-based approach.
Figure \ref{fig:results_convergence_plot_all_unknown} presents the convergence curves for these datasets, and Table \ref{tab:overview_results_unknown} reports the corresponding numerical results.
As before, best-performing results are highlighted in bold and statistically significant ones are underlined.

\begin{table}[h]
	\fontsize{9pt}{14pt}\selectfont
	\centering
    \caption{Overview of the results, unknown instance types}
	\begin{tabular}{c|c|c|c|c|c|c|c|c}
        \hline
        & \multicolumn{4}{c|}{CAS-PFSP-M5T1} & \multicolumn{4}{c}{CAS-PFSP-M5T3} \\
        \hline
        Method & mean & best & std & $\% \Delta$ & mean & best & std & $\% \Delta$ \\
        \hline
        MA default & $\mathbf{\underline{4.59 \times 10^6}}$ & $\mathbf{4.56 \times 10^6}$ & $2.15 \times 10^4$ & -0.27 & $2.39 \times 10^7$ & $2.37 \times 10^7$ & $1.54 \times 10^5$ & 5.73 \\
        MA tuned & $4.60 \times 10^6$ & $4.56 \times 10^6$ & $3.11 \times 10^4$ & 0.04 & $2.27 \times 10^7$ & $2.24 \times 10^7$ & $1.66 \times 10^5$ & 0.31 \\
        MA-DRL & $4.60 \times 10^6$ & $4.56 \times 10^6$  & $2.91 \times 10^4$ & 0 & $\mathbf{\underline{2.26 \times 10^7}}$ & $\mathbf{2.24 \times 10^7}$ & $1.69 \times 10^5$ & 0 \\
        \hline
        & \multicolumn{4}{c|}{CAS-PFSP-M10T1} & \multicolumn{4}{c}{CAS-PFSP-M10T3} \\
        \hline
        Method & mean & best & std & $\% \Delta$ & mean & best & std & $\% \Delta$ \\
        \hline
        MA default & $1.34 \times 10^7$ & $1.33 \times 10^7$ & $3.82 \times 10^4$ & 1.88 & $7.27 \times 10^7$ & $7.14 \times 10^7$ & $7.34 \times 10^5$ & 10.08 \\
        MA tuned & $1.32 \times 10^7$ & $1.31 \times 10^7$ & $5.91 \times 10^4$ & 0.66 & $6.81 \times 10^7$ & $6.74 \times 10^7$ & $3.92 \times 10^5$ & 3.18 \\
        MA-DRL & $\mathbf{\underline{1.31 \times 10^7}}$ & $\mathbf{1.31 \times 10^7}$  & $4.24 \times 10^4$ & 0 & $\mathbf{\underline{6.60 \times 10^7}}$ & $\mathbf{6.54 \times 10^7}$ & $3.90 \times 10^5$ & 0 \\
        \hline
        & \multicolumn{4}{c|}{CAS-PFSP-M15T1} & \multicolumn{4}{c}{CAS-PFSP-M15T3} \\
        \hline
        Method & mean & best & std & $\Delta\%$ & mean & best & std & $\Delta\%$ \\
        \hline
        MA default & $2.02 \times 10^7$ & $2.00 \times 10^7$ & $9.73 \times 10^4$ & 5.29 & $1.02 \times 10^8$ & $1.02 \times 10^8$ & 0 & 4.24 \\
        MA tuned & $1.94 \times 10^7$ & $1.93 \times 10^7$ & $8.04 \times 10^4$ & 1.46 & $1.01 \times 10^8$ & $9.93 \times 10^7$ & $7.16 \times 10^5$ & 2.32 \\
        MA-DRL & $\mathbf{\underline{1.92 \times 10^7}}$ & $\mathbf{1.90 \times 10^7}$  & $6.47 \times 10^4$ & 0 & $\mathbf{\underline{9.83 \times 10^7}}$ & $\mathbf{9.74 \times 10^7}$ & $4.84 \times 10^5$ & 0 \\
        \hline
    \end{tabular}
    \label{tab:overview_results_unknown}
\end{table}

In contrast to the outcomes on the known instance types, the results on the unknown datasets show a clearer advantage for the DRL-based approach.
With the exception of the simplest unknown dataset (M5T1), MA-DRL consistently achieves the best mean and best objective values relative to the tuned baseline.
These improvements are statistically significant, although for the moderately sized datasets their magnitude remains modest.

For example, in M5T3, MA-DRL improves over the tuned configuration by just $0.3 \, \%$, indicating that while the learned policy does transfer beyond the training distribution, the benefit is limited when the problem characteristics deviate only mildly from those seen during training.
As the combinatorial complexity increases, however, the advantage becomes more substantial.
In M10T3, MA-DRL improves upon the tuned baseline by approximately $3.2 \, \%$, and in M15T3 by about $2.3 \, \%$.
These gains are not only statistically significant but also practically meaningful, especially for the larger instances where even small percentage improvements translate into large absolute reductions in objective value.

The contrast with the default configuration is naturally larger, especially for the more complex datasets, with improvements reaching up to $10 \, \%$ in M10T3.
The only dataset in which MA-DRL does not outperform the static methods is M5T1, the simplest of the unknown datasets.
Here, the default configuration slightly outperforms both tuning and DRL, although this difference is extremely small in practical terms.

Taken together with the increasing gains observed for larger and more complex unknown instances, these findings strongly suggest that problem complexity is a central driver of when learning becomes beneficial.
When instances closely resemble those encountered during training, such as in the relatively simple M5T1 dataset, the DRL-based method offers little or no advantage over tuned static parameter settings.
However, as soon as the instances deviate more substantially from the training distribution and their combinatorial complexity increases, the advantage of MA-DRL becomes both clearer and more pronounced.
This notion of difference is therefore defined here in terms of these structural instance characteristics and should be interpreted as an empirical observation rather than a predictive criterion.
Detailed per-instance results are available at \url{https://github.com/ugent-isye/MA-DRL-CAS-PFSP}.

\subsubsection{Analysis of learned parameter control policy}

To better understand the behavior of the learned DRL-based parameter control policy, we analyze the evolution of the selected parameter values throughout the execution of the algorithm.
In particular, we investigate whether the policy continuously adapts the parameters over time or converges to a stable configuration.

\begin{figure}[h]
	\begin{subfigure}{.49\columnwidth}
		\centering
		\includegraphics[clip, trim={0.4cm 0cm 0.4cm 0.5cm}, width=\textwidth]{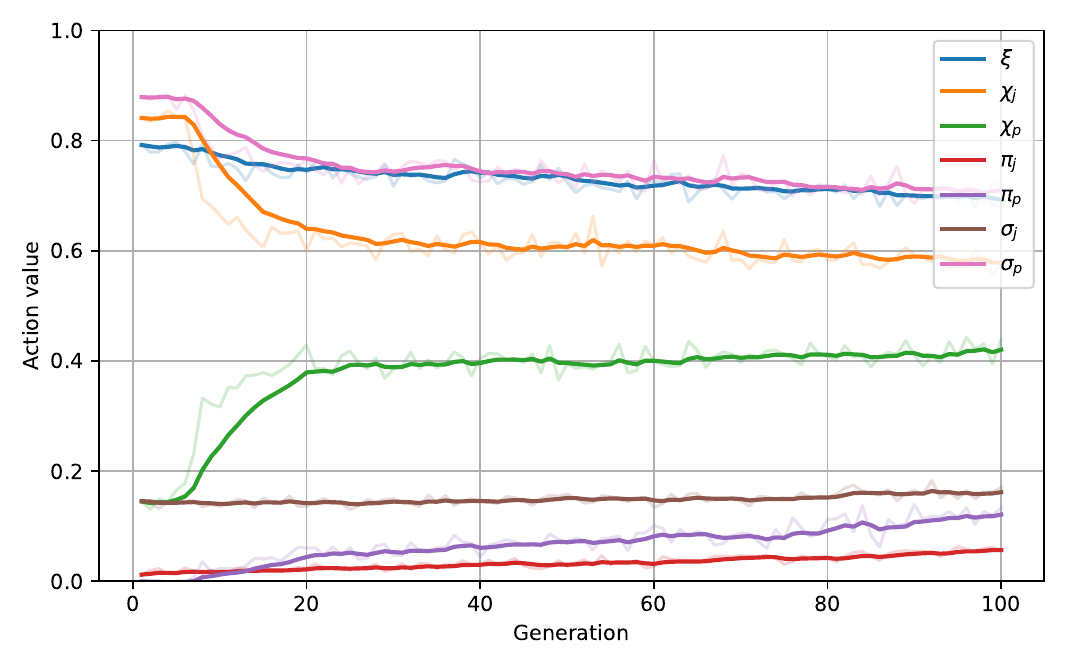}
		\caption{Action values produced by the DRL agent}
		\label{fig:parameters_evolution_rescaled}
	\end{subfigure}
	\begin{subfigure}{.49\columnwidth}
		\centering
		\includegraphics[clip, trim={0.4cm 0cm 0.4cm 0.5cm}, width=\textwidth]{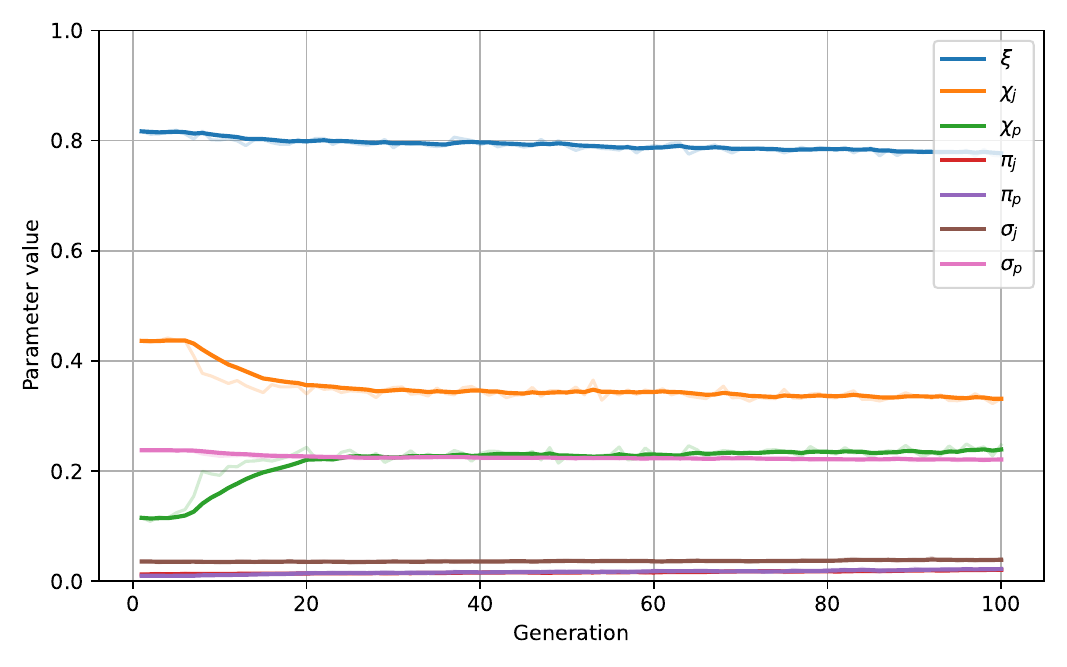}
		\caption{Rescaled parameter values}
		\label{fig:parameters_evolution_true}
	\end{subfigure}
	\caption{Evolution of parameter values during the search}
	\label{fig:parameters_evolution}
\end{figure}

Figure \ref{fig:parameters_evolution} presents the average parameter trajectories over all training instances and runs.
Two complementary representations are shown.
The left plot depicts the action values in the $[0,1]$ space produced by the DRL agent, while the right plot shows the corresponding parameter values after rescaling to their operational ranges.
For each parameter, both the raw trajectories and a smoothed trend (with $\alpha=0.8$) are reported.

The results reveal a clear two-phase behavior of the learned policy.
During the initial stage of the search (approximately the first 20 generations), several parameters are actively adjusted.
In particular, the crossover-related parameters $\chi_j$ and $\chi_p$ exhibit pronounced and systematic changes, indicating that the agent dynamically adapts the exploration--exploitation balance in the early phase of the search.
Other parameters, such as $\xi$ and $\sigma_p$, show more moderate but consistent trends, while parameters like $\pi_j$, $\pi_p$, and $\sigma_j$ remain nearly constant throughout the run.

After this initial adaptation phase, the parameter values stabilize and exhibit only minor variations for the remainder of the search.
This behavior is consistently visible in both representations and indicates that the DRL policy converges to a near-static configuration, rather than continuously adapting the parameters at every generation.

Importantly, the normalized action trajectories confirm that this behavior is not an artifact of the rescaling applied to the action space.
Instead, it reflects the actual decisions produced by the DRL agent: the policy initially explores different regions of the action space and subsequently settles into a relatively stable region.

These observations provide insight into the role of dynamic parameter control in the proposed approach.
They suggest that the DRL agent primarily learns how to adjust parameter values at the beginning of the search, while relying on relatively stable configurations in later stages.
However, further analysis would be required to determine how this behavior varies across instance types and to what extent it contributes to the observed performance differences.

\section{Conclusion}
\label{section:conclusion}

In this paper, we investigated when using learning-based methods for Dynamic Algorithm Configuration is practically beneficial for a real-world combinatorial optimization problem.
Motivated by the high computational cost of classical tuning and the limited generalization capabilities of static configurations and classical dynamic control strategies, we developed a DRL-based framework to dynamically control the parameters of a memetic algorithm for carbon-aware scheduling.
Our aim was not only to test whether DRL can outperform static parameter settings, but explicitly to identify the conditions under which using learning-based DAC is preferable to conventional tuning.

To this end, we constructed ten instance datasets spanning a wide range of combinatorial complexity, and we distinguished between known and unknown instance types.
The former include the relatively simple instances used for both training the DRL agent and tuning the static algorithm, while the latter represent  more complex and previously unseen instance types used exclusively for evaluation.
This separation enabled a controlled analysis of generalization, allowing us to assess whether a policy learned on small, cheap-to-train instances can transfer effectively to larger, more complex scheduling instances, where both tuning and training would be computationally very expensive.

We compared our DRL method to two baselines: (i) a static memetic algorithm with default parameter settings, as commonly used in the literature, and (ii) a static memetic algorithm with tuned parameters.
To ensure a fair comparison, we designed a tuning strategy that mirrors the DRL training process as closely as possible.
In particular, both training and tuning operated on the same set of 12 representative training instances and were allocated the same total computational budget.

From our computational experiments, three key insights emerge.
First, we show that learning-based DAC is not always superior to static algorithm configuration.
Across the known instance types, i.e., those whose characteristics closely match the training and tuning set, the DRL-based approach performs similarly to the static tuned method.
Even when small numeric differences were significant, their magnitudes were negligible, suggesting that learning-based control methods offer no practical benefit when the deployment environment is well represented during training.

Second, we find that the usefulness of learning depends crucially on how different the evaluation instances are from those used for training in terms of their combinatorial complexity.
In this work, such differences are characterized more concretely by structural features such as the number of machines, jobs, and resulting operations, which directly determine the size of the search space.

When instances are similar to the training set, such as in the simplest unknown dataset M5T1, the learned policy provides little or no advantage over tuned static parameters.
However, as soon as problem characteristics begin to diverge, and especially as the number of machines and operations increases, MA-DRL begins to outperform static tuning.
On the more complex unknown instance types (e.g., M10T3 and M15T3), improvements up to more than $3 \, \%$ are observed, showing a clear advantage over static tuning.

We note, however, that this notion of instance similarity is defined in a descriptive rather than predictive manner.
While our results empirically demonstrate how performance varies with increasing problem complexity, they do not provide a formal mechanism to predict a priori when learning-based DAC will be beneficial.

Third, our results show that learned parameter control policies generalize across problem scales.
Despite being trained exclusively on small instances, the DRL-based controller performs robustly on much more complex unseen problems where training would be highly expensive due to the increased runtimes.
This validates the central premise of our study: that effective parameter-control policies can be learned on small instances and subsequently applied to far more complex problems, avoiding the cost of tuning on large-scale instances.

In addition to these main findings, an analysis of the learned parameter control policy provides further insight into how the DRL-based approach operates during the search.
The results suggest that the agent follows a two-phase behavior, with an initial adaptation of parameter values during the early generations, followed by a largely stable configuration in later stages.
While this observation is based on aggregated data and does not distinguish between instance types, it indicates that dynamic parameter control may primarily act through early-stage adaptation rather than continuous adjustments throughout the entire run.

Taken together, these findings highlight both the advantages and limitations of the proposed DRL-based dynamic algorithm configuration approach.
On the one hand, the method enables the learning of parameter control policies that generalize across instance types and problem scales, 
and allows the algorithm to adapt its behavior dynamically during the search.
On the other hand, this added flexibility comes with increased modeling complexity and does not necessarily translate into performance gains in all settings.
In particular, the results show that when instance characteristics are well represented during training, the learned policies perform comparably to well-tuned static configurations, providing no clear advantage but also no disadvantage in terms of solution quality under equal computational budgets.

On a final note, while the observed relative improvements on the largest instances are modest in percentage terms, it is important to interpret these results in the context of large-scale combinatorial optimization, where even small relative gains can translate into substantial absolute improvements in objective value.
More importantly, these gains are achieved in a setting where the control policy is trained exclusively on small instances and transferred to significantly more complex ones, in contrast to classical tuning approaches, which typically require re-optimization for each new instance type or scale.
Nevertheless, while the computational effort required for training is comparable to that of classical tuning under the considered experimental setup, DRL introduces additional modeling complexity and design choices, and the benefits of learning-based DAC are not universal.
In particular, when instance characteristics are well represented during tuning or when problem instances are relatively simple, the additional investment in developing and training the DRL model may not be justified.
The proposed approach is therefore most relevant in settings where problem complexity is high and instance characteristics vary beyond what can be efficiently captured through static tuning.

\subsection{Directions for future work}
Future work could further deepen these findings in several ways.
A first direction is to apply instance space analysis \cite{Smith-Miles2014} to quantify how instance features influence algorithm performance.
In the present work, differences between instances are defined primarily in terms of problem size and combinatorial complexity, which provides a useful but coarse characterization of instance similarity.
Instance space analysis would help clarify why certain methods outperform others under specific conditions and would also support the construction of more effective training sets by ensuring that a sufficiently diverse range of instance characteristics is represented. 

A second promising direction is to investigate whether the generalization capabilities demonstrated here extend beyond a single problem class.
Since the learned policy acts on the search dynamics of the underlying evolutionary algorithm rather than on a specific problem, it is plausible that a policy trained on one class of combinatorial optimization problems could be deployed on another that uses the same EA framework.
Extending this idea further, one could aim to identify a representative set of problem classes on which to train a policy, with the goal of obtaining a parameter control strategy that generalizes across a wide spectrum of combinatorial optimization problems, thereby making DRL-based DAC a reusable, problem-independent module for evolutionary computing.

A third direction is to extend the proposed approach to explicitly handle multi-objective optimization formulations.
In many realistic scheduling settings, multiple conflicting objectives such as makespan, energy cost, and carbon emissions must be considered simultaneously.
Incorporating such trade-offs would allow the proposed framework to be applied in more complex and realistic decision-making contexts, but also introduces additional challenges in terms of solution representation and policy learning.

\section*{Acknowledgments}
This work was funded by the Research Foundation - Flanders (FWO), with grant number: V423124N.

\begingroup
    \setlength{\bibsep}{10pt}
    \setstretch{1}
	\bibliographystyle{apa} 
    \bibliography{refs_CAS}
\endgroup

\end{document}